\documentclass[12pt]{article}
\usepackage{amssymb,amsmath,color}

\def \R{{\hbox{\vrule width 0.6pt height 6.8pt depth -.2pt\kern-0.2pt
R}}}
\def \P{{\hbox{\vrule width 0.6pt height 6.8pt depth -.2pt\kern-0.2pt
P}}}

\def\m {{\eta}}

\def \no { \noindent}
\def \er {\mathbb R}

\def\t0{{t}}

\def \m {\eta}

\def\grad{\nabla}
\def\div{\, \mbox{div}\,  }

\def\MOD#1{{|\kern -.16em |\kern -.16em | #1 | \kern -.16em |\kern
 -.16em |}}
\def \epsilon {\varepsilon}

\def\ds{\displaystyle}
\newtheorem{theo}{\bf THEOREM}
\newtheorem{lem}{\bf LEMMA}[section]
\newtheorem{pro}[lem]{\bf PROPOSITION}

\newtheorem{rem}[lem]{\bf REMARK}

\numberwithin{equation}{section}

\def\B1{B_{1/2}}
\def\Box{\hfill\rule{2.5mm}{2.5mm}}

\def\R{{\mathbb {R}}}

\def\build#1_#2^#3{\mathrel{
\mathop{\kern 0pt#1}\limits_{#2}^{#3}}}

\def\h1{\mathop{\rm H^1_{\rm loc,\rm u}}}

\def\l2{\mathop{\rm L^2_{\rm loc,\rm u}}}

\def\t0{{t}}

\def\Box{\hfill\rule{2.5mm}{2.5mm}}

\begin{document}
\bibliographystyle{siam}

\title{
Blow-up results  for semilinear wave equations in  the
super-conformal case  }

\author{M.A. Hamza\\
{\it \small  {Facult\'e des Sciences de Tunis}}\\
H. Zaag \\
 {\it \small   {CNRS UMR 7539 LAGA  Universit\'e Paris 13}} }

\maketitle

\begin{abstract} We consider
the semilinear wave equation in higher dimensions with power nonlinearity in the super-conformal range, and its perturbations with lower order terms, including the Klein-Gordon equation. We improve the upper bounds on blow-up solutions previously obtained by Killip, Stovall and Vi\c san \cite{KSV}. Our proof uses the similarity variables' setting. We consider the equation in that setting as a perturbation of the conformal case, and we handle the extra terms thanks to the ideas we already developed in \cite{HZJHDE} for perturbations of the pure power case with lower order terms.
\end{abstract}

\vspace{0.5cm}


 \vspace{0.4cm}

\section{\bf Introduction}
This paper is devoted to the study of blow-up solutions for the
following semilinear  wave equation:
\begin{equation}\label{gen}
\left\{
\begin{array}{l}
\partial_t^2 u =\Delta u+|u|^{p-1}u+f(u)+g(x,t,\nabla u,\partial_t u ),\\
\\
(u(x,0),\partial_tu(x,0))=(u_0(x),u_1(x))\in  H^{1}_{loc}(\er^N)\times L^{2}_{loc}(\er^N),\\
\end{array}
\right.
\end{equation}
in spatial dimensions $N\ge 2$, where $u(t):x\in \er^N \rightarrow
u(x,t)\in \er$ and $p_c<p<p_S$, where $ p_c\equiv 1+\frac4{N-1}$ is
the conformal critical exponent and  $p_S\equiv 1+\frac4{N-2}$
is the Sobolev critical exponent. Moreover, we take
 $f:\R\rightarrow \R $ and $g:\R^{2N+2}\rightarrow \R $
${\cal {C}}^1$ functions satisfying
\begin{eqnarray*}
(H_f)& |{f(u)}|\le M(1+|u|^q), \  &{\textrm {for all }}\ u\in \R \  \qquad{{\textrm {with}}}\ \ (q<p,\ \ M>0),\\
(H_g)& |{g(x,t,v,z)}|\le M(1+|v|+|z|), &  {\textrm {for all }}\
x,v\in \R^N, t,z\in \R \  \ {{\textrm {with}}}\  (M>0).
\end{eqnarray*}

\bigskip
\no We would like to mention that equation \eqref{gen}
encompasses the case of the following nonlinear Klein-Gordon
equation
\begin{equation}\label{KG}
 \partial^2_{t}u=\Delta u+ |u|^{p-1} u-u,\qquad  (x,t)\in \er^N\times [0,T).
\end{equation}
In order to keep our analysis clear, we only give the proof  for the following non perturbed equation
\begin{equation}
\label{1} 
 \partial_t^2u=
\Delta u+ |u|^{p-1} u,\qquad  (x,t)\in \er^N\times [0,T), 
\end{equation}
 and refer the reader  to  \cite{HZlyap09} and
\cite{HZJHDE} for straightforward adaptations to
 equation \eqref{gen}.

\bigskip

 \no The Cauchy problem of equation ({\ref{1}}) is solved in
$H^{1}_{loc}\times L^{2}_{loc}$. This follows from the finite speed
of propagation and the the wellposedness in $H^{1} \times L^{2}$,
valid whenever $ 1< p< p_S$. 
The existence of
blow-up solutions for the associated ordinary differential equation
of (\ref{1})  is a classical result. By using the finite speed of
propagation, we conclude that there exists a
 blow-up solution $u(t)$ of (\ref{1}) which depends non trivially on the space variable.
In this paper, we consider  a blow-up solution $u(t)$ of (\ref{1}),
we define (see for example Alinhac \cite{Apndeta95} and
\cite{Afle02})  $\Gamma$ as the graph of a function $x \mapsto T(x)$
such that the domain of definition of  $u$  is given by
$$D_u=\{(x,t)\ \ \big |t<T(x)\}.$$
The set $D_u$ is called the maximal influence domain of $u$.
Moreover, from the finite speed of propagation, $T$ is a
$1$-Lipschitz function. 
The 
 graph
$\Gamma$  is called 
the blow-up graph of $u$. \\
Let us first  introduce the following non-degeneracy condition for
$\Gamma$. If we introduce for all $x\in\er^N$, $t\le T(x)$ and
$\delta>0$, the cone
\begin{equation} \label{cone}
C_{x,t,\delta}=\{(\xi,\tau)\neq(x,t)|0\le \tau\le t-\delta|\xi-x|\},
\end{equation}
then our non degeneracy condition is the following: $x_0$ is a non
characteristic point if
\begin{equation} \label{cone1}
\exists\delta_0=\delta_0(x_0)\in (0,1) \ \textrm{ such that}\  u \
\textrm{ is defined on }\ C_{x_0,T(x_0),\delta_0}.
\end{equation}
\no We aim at studying  the growth estimate  of $u(t)$ near the
space-time blow-up graph in the  super-conformal   case (where
$p_c<p<p_S$).

Let us briefly mention some results concerning the blow-up  rate of
solutions of  semilinear wave equations.
  The first result valid for general solutions is due to Merle and Zaag   in \cite{MZimrn05}
(see also \cite{MZajm03} and \cite{MZma05}) who   proved, that  if
$1<p\le p_c$ and $u$ is a solution of (\ref{1}), then   the growth
estimate near the space-time blow-up graph is given by the
associated ODE. In \cite{HZlyap09} and \cite{HZJHDE}, we extend the
result of Merle and Zaag  to perturbed equations of type (\ref{gen})
under some reasonable growth estimates on $f$ and $g$ in (\ref{gen})
(see hypothesis $(H_f)$ and $(H_g)$).
 Note that, in all  these papers,
 the method  crucially relies on the existence of a Lyapunov
functional in similarity variables established by Antonini and Merle
\cite{AMimrn01}.
 Recently, Killip, Stovall and Vi\c san in \cite{KSV} have
shown, among other results, that the results of
 Merle and Zaag  remain valid for the semilinear Klein-Gordon equation (\ref{KG}).
 Moreover, they  consider also the case where $p_c<p<p_S$ and prove
 that, if $u $   is a solution of ({\ref{KG}}),  then for all  $x_0\in \er^N$, there
exists $K>0$ such that, for all $t\in [0,T(x_0))$,
\begin{eqnarray}\label{v1}
(T(x_0)-t)^{-\frac{(p-1)N}{p+3}}\int_{B(x_0,\frac{T(x_0)-t}2)} \! u^2(x,t)\mathrm{d}x &\le&
K,
\end{eqnarray}
and for all   $t\in (0,T(x_0)]$,
\begin{eqnarray}\label{v2}
\int_{T(x_0)-t}^{T(x_0)-\frac{t}2}\int_{B(x_0,\frac{T(x_0)-\tau }2)}
\!\Big(|\grad u(x,\tau )|^{2}+|\partial_tu(x,\tau
)|^2\Big)\mathrm{d}x\mathrm{d}\tau&\le&K.
\end{eqnarray}
Moreover, if   $x_0$ is a non characteristic point,  then they use
 a covering argument to obtain  the same
estimates with   the ball $B(x_0,\frac{T(x_0)-\tau}2)$ replaced by the ball
$B(x_0,T(x_0)-\tau )$ in the inequalities (\ref{v1}) and (\ref{v2}).

\bigskip

\no Here, we obtain a better result
 thanks to a
different method based  on the use of  self-similar variables. This
method allows us to improve the results of \cite{KSV} as we state in the following:
\begin{theo}\label{t0}\textbf{(Growth estimate near the blow-up surface for Eq. ({\ref{gen}}))}.\\
If $u $   is a solution of ({\ref{gen}}) with blow-up graph
$\Gamma:\{x\mapsto T(x)\}$, then for all  $x_0\in \er^N$  and  $t\in
[0,T(x_0))$, we have
\begin{eqnarray}\label{t1}
(T(x_0)-t)^{\frac{-(p-1)N}{p+3}}\int_{B(x_0,T(x_0)-t)} \!
u^2(x,t)\mathrm{d}x \ \ \rightarrow 0 ,\ \ \textrm{ as}\ \ \
t\rightarrow T(x_0).
\end{eqnarray}
Moreover, for all   $t\in (0,T(x_0)]$, we have
\begin{eqnarray}\label{t2}
\int_{T(x_0)-t}^{T(x_0)-\frac{t}2}\int_{B(x_0,\frac{T(x_0)-\tau }2)}
\!|\partial_tu(x,\tau )|^2\mathrm{d}x \mathrm{d}\tau&\le&K_1,
\end{eqnarray}and
\begin{eqnarray}\label{t22}
\int_{T(x_0)-t}^{T(x_0)-\frac{t}2}\int_{B(x_0,\frac{T(x_0)-\tau }2)}
\!|\grad u(x,\tau )|^{2}\mathrm{d}x \mathrm{d}\tau &\le&K_1.
\end{eqnarray}
If in addition  $x_0$ is a non characteristic point, then we have
for all $t\in (0,T(x_0)]$,
\begin{eqnarray}\label{t4}
\int_{T(x_0)-t}^{T(x_0)-\frac{t}2}\int_{B(x_0,{T(x_0)-\tau })}
\!\Big(|\grad u(x,\tau )|^{2}+|\partial_tu(x,\tau
)|^2\Big)\mathrm{d}x\mathrm{d}\tau\ \ \rightarrow 0, \ \textrm{ as}\
\ t\rightarrow 0.\qquad
\end{eqnarray}
Moreover, we have
\begin{eqnarray}\label{limener}
&&\frac{T(x_0)-t}2\int_{B(x_0,{T(x_0)-t
})} \!\Big( |\grad u(x,t )|^{2}- \big(\frac{x-x_0}{T(x_0)-t}.\grad u(x,t )\big)^2\nonumber\\
&&+|\partial_tu(x,t )|^2-\frac1{p+1} \!|
u(x,t)|^{p+1}\Big)\mathrm{d}x\rightarrow 0, \ \textrm{ as}\ \
t\rightarrow T(x_0).
\end{eqnarray}
\end{theo}
\begin{rem}{}
i)Let us remark that,  we have the following   lower bound which
follows  from standard techniques (scaling arguments, the
wellposedness in $H^1(\er^N )\times L^2(\er^N)$, the finite speed of
propagation and the fact that $x_0$ is a non characteristic point):
there exist
 $\varepsilon_0>0$, such that
\begin{eqnarray*}
&&0<\varepsilon_0\le (T(x_0)-t)^{\frac{2}{p-1}}\frac{\|u(t)\|_{L^2(B(x_0,{T(x_0)-t}))}}{ (T(x_0)-t)^{\frac{N}{2}}}\nonumber\\
&&+ (T(x_0)-t)^{\frac{2}{p-1}+1}\Big
(\frac{\|\partial_tu(t)\|_{L^2(B(x_0,{T(x_0)-t}))}}{
(T(x_0)-t)^{\frac{N}{2}}}+
 \frac{\|\grad u(t)\|_{L^2(B(x_0,{T(x_0)-t}))}}{ (T(x_0)-t)^{\frac{N}{2}}}\Big
 ).
\end{eqnarray*}
ii) In  Theorem 1, we improve recent results of
 Killip, Stovall and Vi\c san in \cite{KSV}. More precisely,
  we obtain a better estimate in (\ref{t1}) and  if $x_0$ is non
characteristic point we have the better estimate  (\ref{t4}).\\
iii) Up to a time dependent factor, the expression in (\ref{limener})
is equal to the main terms of the energy  in similarity variables
(see (\ref{E3})).
However, even with this improvement, we think that our estimates are still not optimal. \\
iv) The constant $K_1$, and the rate of convergence to $0$ of the
different quantities in the previous theorem  and in the whole
paper, depend only on $N$, $p$ and the upper bound on $T(x_0)$,
$1/{T(x_0)}$, and the initial data $(u_0,u_1)$ in $
H^{1}(B(x_0,2T(x_0)))\times L^{2}(B(x_0,2T(x_0)))$, together with
$\delta_0(x_0)$ if $x_0$ is non characteristic point.

\end{rem}

\bigskip
Our method relies on the estimates in similarity variables
introduced in \cite{AMimrn01} and used in \cite{MZajm03},
\cite{MZimrn05} and \cite{MZma05}. More precisely, given $(x_0,T_0)$
such that $0< T_0\le T(x_0)$, we introduce  the following
self-similar change of variables:
\begin{equation}\label{scaling}
y=\frac{x-x_0}{T_0-t},\qquad s=-\log (T_0-t),\qquad
u(x,t)=\frac{1}{(T_0-t) ^{\frac{2}{p-1}}}w_{x_0,T_0}(y,s).
\end{equation}
 This change of variables
transforms the backward light cone with vortex $(x_0,T_0)$ into the
infinite cylinder $(y,s)\in B\times [-\log T_0,+\infty)$. In the new
set of variables $(y,s),$ the behavior of $u$ as $t \rightarrow T_0$
is equivalent to the behavior of $w$ as $s \rightarrow +\infty$.

\bigskip

From (\ref{1}), the  function $w_{x_0,T_0}$  (we write $w$ for
simplicity) satisfies the following equation for all $y\in B\equiv
B(0,1)$ and $s\ge -\log T_0$:
\begin{eqnarray}\label{A}
\partial_{s}^2w+\frac{p+3}{p-1}\partial_sw+2y.\grad \partial_sw&=&
\sum_{i,j}(\delta_{i,j}-y_iy_j) \partial^2_{y_i,y_j} w-\frac{2p+2}{(p-1)^2}y.\grad
w\nonumber\\
&&-\frac{2p+2}{(p-1)^2}w+|w|^{p-1}w.
\end{eqnarray}
Putting this equation in the following form
\begin{eqnarray}\label{C}
\partial_{s}^2w&=&\div( \grad w-(y.\grad w)y)+2\eta
 y.\grad w
-\frac{2p+2}{(p-1)^2}w+|w|^{p-1}w\\
&&-\frac{p+3}{p-1}\partial_s w-2y.\grad \partial_sw,\
  \forall y\in B\ {\textrm{ and}} \ s\ge -\log T_0,\nonumber
\end{eqnarray}
where
\begin{equation}\label{eta}
\eta=\frac{N-1}{2}-\frac2{p-1}=\frac2{p_c-1}-\frac2{p-1}>0,
\end{equation}
 the key idea of our paper is to view this equation as a perturbation of the conformal case (corresponding to $\eta=0$) already treated in \cite{HZJHDE} with the term $2\eta y \cdot \nabla w$. Of course, this term is not a lower order term with respect to the nonlinearity. For that reason, we will have exponential growth rates in the $w$ setting. Let us emphasize the fact that our analysis is not just a trivial adpatation of our previous work \cite{HZJHDE}.

\bigskip

\no The equation (\ref{A}) will be studied in the Hilbert  space $\cal H$
$${\cal H}=\Big \{(w_1,w_2), |
\displaystyle\int_{B}\Big ( w_2^2 +|\grad w_1|^2(1-|y|^2)+w_1^2\Big
) {\mathrm{d}}y<+\infty \Big \}.$$
%
%
%
%
%
In the conformal case where $p=p_c$, Merle and Zaag \cite{MZimrn05}
proved that
\begin{eqnarray}
\label{f} E_0(w)=\displaystyle\int_{B}\Big (
\frac{1}{2}(\partial_sw)^2 +\frac{1}{2}|\grad
w|^2-\frac{1}{2}(y.\grad w)^2+\frac{p+1}{(p-1)^2}w^2
-\frac{|w|^{p+1}}{p+1}\Big ) {\mathrm{d}}y,\quad
\end{eqnarray}
 is a Lyapunov functional for equation (\ref{A}).
 When $p>p_c$, we introduce
\begin{eqnarray}\label{E1}
E(w)&=&\displaystyle {E_0(w)+I(w),}
\end{eqnarray}
where
\begin{eqnarray}
\label{E2}
 I(w)&=& -\eta\displaystyle\int_{B}w\partial_s w {\mathrm{d}}y+ \frac{\m N}2 \displaystyle\int_{B}w^2
 {\mathrm{d}}y,
\end{eqnarray}
and $\eta$ is defined in \eqref{eta}.
 Finally, we define the energy function as
\begin{eqnarray}
\label{E3}
 F(w,s)&=&\displaystyle {E(w)e^{-2\eta s}.}
\end{eqnarray}

The  proof of Theorem 1 crucially relies on the fact that  $F(w,s)$ is
 a Lyapunov functional for equation (\ref{A})  on the one hand, and on the other hand, on   a blow-up criterion involving
 $F(w,s)$.
 Indeed, with the  functional $F(w,s)$ and some
more work, we are able to adapt the analysis performed in
\cite{MZimrn05}. In the following, we show that $F(w,s)$  is a
Lyapunov functional:

\begin{pro}\textbf{(Existence  of  a decreasing  functional for Eq. ({\ref{A}}))}.\\
\label{pro} For all $s_2>s_1\ge -\log T_0=s_0$, the functional
$F(w,s)$ defined in (\ref{E3}) satisfies
\begin{eqnarray}\label{LE}
F(w(s_2),s_2)-F(w(s_1),s_1) &=&
- \int_{s_1}^{s_2} e^{-2\m s}\!\int_{\partial B}\!\Big(\partial_s w-\m w\Big)^2{\mathrm{d}}\sigma {\mathrm{d}}s\nonumber\\
&&-  \frac{\m (p-1)}{p+1}\int_{s_1}^{s_2} e^{-2\m
s}\int_{B}|w|^{p+1}{\mathrm{d}}y{\mathrm{d}}s.
\end{eqnarray}
Moreover, for all $s\ge s_0$, we have  $F(w,s)\ge0$.
\end{pro}


\bigskip

This paper is organized as follows:  In section 2, we prove
Proposition \ref{pro}. Using this result, we prove Theorem \ref{t0}
in section 3.

\section{Existence  of  a decreasing  functional for equation (\ref{A}) and a
blow-up criterion} Consider $u $ a solution of ({\ref{1}}) with
blow-up graph $\Gamma:\{x\mapsto T(x)\}$, and consider its
self-similar transformation $w_{x_0,T_0}$ defined at some scaling
point $(x_0,T_0)$ by (\ref{scaling}) where $T_0\le T(x_0)$.
 This section is devoted to the proof of
Proposition \ref{pro}. We proceed in two parts:
\begin{itemize}
\item  In subsection 2.1, we show the existence  of  a decreasing  functional for equation (\ref{A}).
\item In subsection 2.2, we prove a blow-up criterion involving this functional.
\end{itemize}
\subsection{Existence  of  a decreasing  functional for equation
(\ref{A})} In this subsection, we prove  that the functional $F(w,s)$
defined in (\ref{E3})  is decreasing. More precisely we prove that
the functional $F(w,s)$ satisfies the inequality (\ref{LE}).
%
%
%
%
%
Now we state two lemmas which are crucial for the proof. We
begin with   bounding the time derivative of $ E_0(w)$ defined in
(\ref{f}) in the following lemma.
\begin{lem}
\label{energylyap0} For all $s\ge -\log T_0$, we have
\begin{eqnarray}\label{energylem21}
 \frac{d}{ds}(E_0(w))&=& - \int_{\partial B} (\partial_s w)^2
{\mathrm{d}}\sigma + 2 \m \int_{B} (\partial_s w)^2 {\mathrm{d}}y
+2\m \int_{B} \partial_s w(y.\grad w){\mathrm{d}}y.
\end{eqnarray}
\end{lem}
{\it{Proof}}:  Multiplying $(\ref{C})$ by $\partial_s w $ and
integrating over the ball $B$, we obtain for all $s\ge -\log T_0$,
\begin{eqnarray}\label{energyE0}
\frac{d}{ds}(E_0(w))&=& -2 \int_{B} \partial_s w(y.\grad
\partial_s w){\mathrm{d}}y
-\frac{p+3}{p-1} \int_{B}(\partial_s w)^2{\mathrm{d}}y +2\m \int_{B}
\partial_s w(y.\grad w){\mathrm{d}}y.\nonumber
\end{eqnarray}
Since we see from  integration by parts that
\begin{eqnarray*}
 - 2\int_{ B} \partial_s w(y.\grad \partial_s w)  {\mathrm{d}}y&=& - \int_{ B}y.\grad (\partial_s w)^2  {\mathrm{d}}y=
 N
\int_{B} (\partial_s w)^2 {\mathrm{d}}y - \int_{\partial B}
(\partial_s w)^2{\mathrm{d}}\sigma,
\end{eqnarray*}
this concludes the proof of Lemma
 \ref{energylyap0}.
\Box

\vspace{0.3cm}

\no We are now going to prove the following estimate for the
functional $I(w)$:
\begin{lem}
\label{energylyap1} For all $s\ge -\log T_0$, we have
\begin{eqnarray}\label{energytheta1}
\frac{d}{ds}I(w)&=& 2\m E(w) -2\eta \int_{B}(\partial_s
w)^2{\mathrm{d}}y - \frac{\m
(p-1)}{p+1}\int_{B}|w|^{p+1}{\mathrm{d}}y\nonumber\\
&&-2\m \int_{B}
\partial_sw(y.\grad w) {\mathrm{d}}y  -\m^2 \int_{\partial B} w^2
{\mathrm{d}}\sigma+2\m \int_{\partial B} w\partial_s w
{\mathrm{d}}\sigma.
\end{eqnarray}
\end{lem}
{\it Proof:} Note that  $I(w)$ is a differentiable  function for all
 $s\ge   -\log T_0$ and that
\begin{eqnarray*}
\frac{d}{ds}I(w)&=&- \m\int_{B}(\partial_sw)^2{\mathrm{d}}y-
\m\int_{B}w\partial^2_{s}w{\mathrm{d}}y +\m N \int_{B}w\partial_s w
{\mathrm{d}}y.
\end{eqnarray*}
By using  equation (\ref{C}) and integrating by parts, we have
\begin{eqnarray*}
\frac{d}{ds}I(w)&=& -\eta \int_{B}(\partial_s w)^2{\mathrm{d}}y
+ \m\int_{B}(|\grad w|^2-(y.\grad w)^2){\mathrm{d}}y- \m\int_{B}|w|^{p+1}{\mathrm{d}}y\nonumber\\
&&-2\m^2 \int_{B} w(y.\grad w){\mathrm{d}}y +\m \frac{2p+2}{(p-1)^2}\int_{B}w^2{\mathrm{d}}y\\
&&+2\m \int_{B} w(y.\grad \partial_sw) {\mathrm{d}}y +\m
(\frac{p+3}{p-1}+N) \int_{B}w\partial_sw{\mathrm{d}}y.\nonumber
\end{eqnarray*}
Then by integrating by parts, we have
\begin{eqnarray}\label{e8}
\frac{d}{ds}I(w)&=& -\eta \int_{B}(\partial_s w)^2{\mathrm{d}}y
+ \m\int_{B}(|\grad w|^2-(y.\grad w)^2){\mathrm{d}}y- \m\int_{B}|w|^{p+1}{\mathrm{d}}y\nonumber\\
&&-\m^2 \int_{\partial B}  w^2{\mathrm{d}}\sigma +(\m^2N  +\m \frac{2p+2}{(p-1)^2})\int_{B}w^2{\mathrm{d}}y\\
&&+2\m \int_{\partial B} w \partial_sw {\mathrm{d}}\sigma -2\m
\int_{B} (y.\grad w)
\partial_sw {\mathrm{d}}y  + \m
(\frac{p+3}{p-1}-N) \int_{B}w\partial_sw{\mathrm{d}}y.\nonumber
\end{eqnarray}
By combining (\ref{f}), (\ref{E1}),  (\ref{E2}), (\ref{eta}) and
(\ref{e8}), we conclude the proof of Lemma \ref{energylyap1}.
\Box

\bigskip
\no From Lemmas \ref{energylyap0} and \ref{energylyap1}, we are in a
position to prove   the first part of Proposition \ref{pro}.

\vspace{0.2cm} \no {\it Proof of the first part of  Proposition
\ref{pro}:} From Lemmas \ref{energylyap0} and \ref{energylyap1}, we
obtain for all $s\ge -\log T_0$,
\begin{eqnarray*}
\frac{d}{ds}E(w)&=&2\m E(w)- \int_{\partial B} \Big(\partial_s w-\m
w\Big)^2 {\mathrm{d}}\sigma - \frac{\m
(p-1)}{p+1}\int_{B}|w|^{p+1}{\mathrm{d}}y.
\end{eqnarray*}
Therefore, using  the definition of the functional $F(w,s)$ in
(\ref{E3}), we  write
\begin{eqnarray}\label{m}
\frac{d}{ds}F(w,s)&=&-e^{-2\m s} \int_{\partial B} \Big(\partial_s
w-\m w\Big)^2 {\mathrm{d}}\sigma - \frac{\m (p-1)}{p+1}e^{-2\m s}
\int_{B}|w|^{p+1}{\mathrm{d}}y  .
\end{eqnarray}
By integration, we get (\ref{LE}). This concludes the first part of
the proof of Proposition \ref{pro}. \Box

\subsection{A  blow-up criterion}
\no We finish  the proof  of  Proposition \ref{pro} here. More
precisely, for all $ x_0\in \er^N$ and $ T_0\in (0,T(x_0)]$, we
prove that
\begin{equation}\label{254}
\forall \ s\ge -\log T_0,\ \ \ F(w_{x_0,T_0}(s),s)\ge 0.
\end{equation}
 We  give the proof only  in the case where $x_0$ is a non
characteristic point. Note that the case where $x_0$ is a
 characteristic point  can be  done  exactly as in Appendix A page 119 in \cite{MZjfa07}.

\bigskip

\no

{\it{Proof of the last point of Proposition 1.2}}:  The argument is the
same as in the corresponding part in \cite{AMimrn01}. We write the
proof for completeness. Arguing by contradiction, we assume that
there exists  a non characteristic point $ x_0\in \er^N$, $ T_0\in
(0,T(x_0)]$  and  $ s_1\ge -\log T_0$ such that $F(w(s_1),s_1)<0$, where
$w=w_{x_0,T_0}$. Since the energy $F(w(s),s)$ decreases in time, we
have $F(w(1+s_1),1+s_1)<0$.

\no Consider now for $\delta>0$ the function
$\widetilde{w}^{\delta}(y,s)=w_{x_0,T_0-\delta }(y,s)$. From
(\ref{scaling}), we see that for all $(y,s)\in B\times
[1+s_1,+\infty)$
\begin{equation*}
 \widetilde{w}^{\delta}(y,s)=\frac{1}{(1+\delta e^s)^{\frac2{p-1}}}
w( \frac{y}{1+\delta e^s},-\log(\delta+e^{-s})).
\end{equation*}
\begin{itemize}
\item  (A) Note that $\widetilde{w}^{\delta}$ is defined in  $ B\times [1+s_1,+\infty)$,
whenever $\delta>0$ is small enough such that
$-\log(\delta+e^{-1-s_1})\ge s_1.$
\item (B) By construction, $\widetilde{w}^{\delta}$
is also a solution of equation (\ref{A}).
\item (C) For $\delta$ small enough, we have
$F(\widetilde{w}^{\delta}(1+s_1),1+s_1)<0$ by continuity of the function
$\delta \mapsto F(\widetilde{w}^{\delta}(1+s_1),1+s_1)$.
\end{itemize}
Now, we fix $\delta=\delta_0>0$ such that (A), (B) and (C) hold. Since
 $F(\widetilde{w}^{\delta_0},s)$ is decreasing  in time, we have
 \begin{equation}
\label{255}
\liminf_{s\rightarrow +\infty}F(\widetilde{w}^{\delta_0}(s),s)\le F(\widetilde{w}^{\delta_0}(1+s_1),1+s_1)<0.
 \end{equation}
 Let
us note that we have
\begin{eqnarray}\label{c1}
-\m\int_{B}\widetilde{w}^{\delta_0}\partial_s\widetilde{w}^{\delta_0}
{\mathrm{d}}y&\ge&-\frac{1}{2}
\int_{B}(\partial_s\widetilde{w}^{\delta_0})^2{\mathrm{d}}y-\frac{\m^2}2
\int_{B}(\widetilde{w}^{\delta_0})^2{\mathrm{d}}y
\end{eqnarray}
By (\ref{f}), (\ref{E1}), (\ref{E2}), (\ref{c1}) and the fact that $\m \in [0,1]$,  we deduce
\begin{eqnarray}
E(\widetilde{w}^{\delta_0}(s))&\ge& (\frac{\m N}2-
\frac{\m^2}2)\int_{B}(\widetilde{w}^{\delta_0})^2
{\mathrm{d}}y-\frac{1}{p+1}\int_{B}|\widetilde{w}^{\delta_0}|^{p+1}
{\mathrm{d}}y\nonumber\\
&\ge& -\frac{1}{p+1}\int_{B}|\widetilde{w}^{\delta_0}|^{p+1}
{\mathrm{d}}y.
\end{eqnarray}
So, by (\ref{E3}), we have
\begin{eqnarray}
F(\widetilde{w}^{\delta_0}(s),s) &\ge& -\frac{e^{-2\m
s}}{p+1}\int_{B}|\widetilde{w}^{\delta_0}|^{p+1} {\mathrm{d}}y.
\end{eqnarray}
After a change of variables, we find that
\[
F(\widetilde{w}^{\delta_0}(s),s) \ge -\frac{e^{-2\m s}}
{(p+1)(1+\delta_0e^s)^{\frac{4}{p-1}+2-N}}\displaystyle\int_{B}|w(z,-\log
(\delta_0+e^{-s}))|^{p+1} {\mathrm{d}}z.
\]
Since we have $-\log (\delta_0+e^{-s})\rightarrow -\log \delta_0$
as $s\rightarrow +\infty$ and since $\|w(s)\|_{L^{p+1}(B)}$ is
locally bounded from the fact that $w=w_{x_0,T_0}$ and $x_0$ is  non
characteristic point, by a continuity argument, it follows that the
former integral remains bounded and
\begin{eqnarray}\label{ggg}
F(\widetilde{w}^{\delta_0}(s),s)&\ge& -\frac{Ce^{-2\m
s}}{(1+\delta_0e^s)^{\frac{4}{p-1}+2-N}}\rightarrow 0,
\end{eqnarray}
as $s\rightarrow +\infty$ (use the fact that $\frac{4}{p-1}+2-N-2\m
=1$ and $\m >0$). So, from (\ref{ggg}), it follows that
\begin{eqnarray}\label{d2}
\liminf_{s\rightarrow +\infty}F(\widetilde{w}^{\delta_0}(s),s)\ge 0 .
\end{eqnarray}
\no From (\ref{255}), this is a contradiction. Thus (\ref{254}) holds.
 This concludes  the  proof of Proposition \ref{pro}.
\Box

\section{Proof of Theorem \ref{t0}}
\no Consider $u $ a solution of ({\ref{1}}) with blow-up graph
$\Gamma:\{x\mapsto T(x)\}$. Translating Theorem \ref{t0} in the self-similar setting   
$w_{x_0,T_0}$ (we write $w$ for simplicity) defined
by (\ref{scaling}), our goal becomes the following Proposition:
\begin{pro}
\label{final}
If $u $   is a solution of ({\ref{gen}}) with blow-up graph
$\Gamma:\{x\mapsto T(x)\}$, then for all  $x_0\in \er^N$  and  $T_0\le
T(x_0)$, we have
for all $s\ge s_0=-\log T_0$,
\begin{eqnarray}\label{cor1}
e^{-2\eta s}\int_{s}^{s+1} \!\int_{B}\!\Big((\partial_s
w(y,\tau))^2+|\grad
w(y,\tau)|^2(1-|y|^2)\Big){\mathrm{d}}y{\mathrm{d}}\tau\le K.
\end{eqnarray}
Moreover,
\begin{eqnarray}\label{cor2}
e^{-2\m s}\int_{ B}\!|w(y,s)|^{\frac{p+3}2}{\mathrm{d}}y \ \
\rightarrow 0 ,\ \ \textrm{ as}\ \ \ s\rightarrow +\infty ,
\end{eqnarray}
\begin{eqnarray}\label{cor3}
e^{\frac{-8\m s}{p+3}}\int_{ B}\!|w(y,s)|^2{\mathrm{d}}y \ \
\rightarrow 0 ,\ \ \textrm{ as}\ \ \ s\rightarrow +\infty .
\end{eqnarray}
If in addition  $x_0$ is a non characteristic point, then we have,
\begin{eqnarray}
e^{-2\eta s}\int_{s}^{s+1} \!\int_{B}\!\Big((\partial_s
w(y,\tau))^2
\Big){\mathrm{d}}y{\mathrm{d}}\tau\ \ \ &\rightarrow& 0\label{03bis} ,\\
e^{-2\eta s}\int_{s}^{s+1} \!\int_{B}\!\Big(
|\grad
w(y,\tau)|^2\Big){\mathrm{d}}y{\mathrm{d}}\tau\ \ \ &\rightarrow& 0 ,\label{03}
\end{eqnarray}
as $s\rightarrow +\infty$. Moreover, we have
\begin{eqnarray}\label{04}
F(w,s)\ \rightarrow 0 ,\ \ \textrm{ as}\ \ \ s\rightarrow +\infty .
\end{eqnarray}
\end{pro}
In this section, we prove  Proposition \ref{final} which directly implies Theorem \ref{t0},  as in  the proof of Theorem 1.1, (page 1145)  in \cite{MZimrn05}.

\bigskip

Let us first use Proposition \ref{pro} and the averaging technique of \cite{MZma05} and \cite{MZimrn05} to get the following bounds:
\begin{lem}\label{lem1}
For  all $s\ge s_0=-\log T_0$, we have
\begin{equation}\label{F1}
0\le F(w(s),s)\le F(w(s_0),s_0),
\end{equation}
\begin{eqnarray}\label{F2}
\int_{s_0}^{\infty}e^{-2\m s}
\!\!\int_{B}\!|w(y,s)|^{p+1}{\mathrm{d}}y{\mathrm{d}}s &\le&
\frac{p+1}{\m (p-1)} F(w(s_0),s_0),
\end{eqnarray}
\begin{eqnarray}\label{F3}
\int_{s_0}^{\infty}e^{-2\m s} \!\!\int_{\partial B}\!\Big(\partial
_sw(\sigma ,s)-\m w(\sigma ,s)\Big)^{2}{\mathrm{d}}\sigma
{\mathrm{d}}s &\le& F(w(s_0),s_0).
\end{eqnarray}
If in addition $x_0$ is non characteristic (with a slope $\delta_0\in (0,1)$), then
\begin{eqnarray}\label{cor00001}
e^{-2\m s}\int_{s}^{s+1}\int_{B}\Big (\partial_s w_{x_0,T_0}(y,\tau
)-\lambda(\tau ,s) w_{x_0,T_0}(y,\tau )\Big )^2dy d\tau \rightarrow
0, \ \ {\textrm {as}}\  s\rightarrow +\infty,\ \
\end{eqnarray}
where $0\le \lambda(\tau ,s) \le C(\delta_0)$, for all $\tau \in
[s,s+1]$.
\end{lem}
\no {\it{Proof:}}
  The first three estimates are  a direct
consequence of Proposition \ref{pro}. As for the last estimate, by introducing $f(y,s)=e^{-\m s}w(y,s)$, we see that  the
dispersion estimate (\ref{F3}) can be written as follows:
\begin{eqnarray}\label{F4}
\int_{s_0}^{\infty} \!\!\int_{\partial B}\!\Big(\partial _sf(\sigma
,s)\Big)^{2}{\mathrm{d}}\sigma {\mathrm{d}}s &\le& F(w(s_0),s_0).
\end{eqnarray}
In particular, we have
\begin{eqnarray}\label{F5}
\int_{s}^{s+1} \!\!\int_{\partial B}\!\Big(\partial _sf(\sigma
,\tau)\Big)^{2}{\mathrm{d}}\sigma {\mathrm{d}}\tau \rightarrow 0, \
\ {\textrm {as}}\ \ s\rightarrow +\infty.
\end{eqnarray}
 By exploiting (\ref{F5}) where the space integration is done over the unit sphere, one can use the averaging technique of Proposition 4.2 (page 1147) in   \cite{MZimrn05} to get the same estimate with the space variable integrated over the whole ball $B$.
\Box

\bigskip

From Lemma \ref{lem1}, we are in a position to prove
Proposition \ref{final}

\bigskip

\no {\it{Proof of
Proposition \ref{final}:}}\\
- {\it Proof of \eqref{cor1}}:
 By integrating  the functional $F(w,s)$  defined in  (\ref{E3}) in time between $s$ and $s+1$, we obtain:
\begin{eqnarray}\label{et}
&&\displaystyle\int_{s}^{s+1}\!\!e^{-2\m
\tau}\int_{B}\!\!\Big((\partial_sw)^2+|\grad w|^2(1-|y|^2)\Big)
{\mathrm{d}}y{\mathrm{d}}\tau=  -\frac{2(p+1)}{(p-1)^2}
\displaystyle\int_{s}^{s+1}\!\!e^{-2\m \tau}\int_{B}\!\! w^2
{\mathrm{d}}y{\mathrm{d}}\tau\nonumber \\
&&2\int_{s}^{s+1}\!\!F(w(\tau ),\tau )
d\tau-\displaystyle\int_{s}^{s+1}\!e^{-2\m \tau}\int_{B}\!\!\Big
(|\grad w|^2|y|^2-(y.\grad w)^2\Big )
{\mathrm{d}}y{\mathrm{d}}\tau\\
&&+\frac{2}{p+1}\displaystyle\int_{s}^{s+1}e^{-2\m
\tau}\int_{B}\!\!|w|^{p+1} {\mathrm{d}}y{\mathrm{d}}\tau
 +\underbrace{2\m \int_{s}^{s+1}\!e^{-2\m
\tau}\displaystyle\int_{B}\!\!\big(w\partial_s w
 -\frac{
N}2 w^2\big) {\mathrm{d}}y{\mathrm{d}}\tau }_{A(s)}.\nonumber
\end{eqnarray}
Now, we    control all the terms on the right-hand side of the
relation (\ref{et}):

\no Note that the first term is negative, while the second term is
bounded because of the  bound (\ref{F1}) on the energy $F(w,s)$.
 Since $|y.\grad
w|\le |y| |\grad w|$, we can say that the third is also negative.
Remark that  (\ref{F2}) implies that the fourth term is also bounded.
Finally, it remains only to control the  term  $A(s)$.

\medskip

\no Combining  the Cauchy-Schwarz inequality, the inequality  $ab\le
 \varepsilon a^2+ \frac1{4\varepsilon}b^2$, and the fact that $N\ge 2$ and $\m\in [0,1]$, we write
\begin{eqnarray}\label{et1}
 A(s)&\le& \frac{1}2 \int_{s}^{s+1}\!\!e^{-2\m
\tau}\displaystyle\int_{B}\!\!(\partial_s w)^2
{\mathrm{d}}y{\mathrm{d}}\tau.
\end{eqnarray}
Now, we are able to conclude the proof of the inequality
(\ref{cor1}). For this, we combine (\ref{et}), (\ref{et1}) and  the above-mentioned arguments for
the first four terms 
to get
\begin{eqnarray}\label{et4}
\int_{s}^{s+1}\!\!\!e^{-2\m \tau}\!\!\int_{B}\!\!\!\Big (
(\partial_sw)^2+|\grad w|^2(1-|y|^2)\Big)
{\mathrm{d}}y{\mathrm{d}}\tau \le K + \frac{1}{2}
\!\int_{s}^{s+1}\!\!\!e^{-2\m \tau}\!\!\int_{B}\!\!\!(\partial_s
w)^2 {\mathrm{d}}y{\mathrm{d}}\tau.\quad
\end{eqnarray}
The desired bound in (\ref{cor1}) follows then from (\ref{et4}).

\medskip

\no - {\it Proof of \eqref{cor2}}:
 Using the mean value theorem, we derive the existence of $\sigma(s)\in [s,s+1]$ such that
\begin{eqnarray}\label{s1}
\int_{ B}\!|w(y,\sigma (s))|^{\frac{p+3}2}{\mathrm{d}}y
&=&\int_{s}^{s+1}\int_{ B}\!|w(y,\tau
)|^{\frac{p+3}2}{\mathrm{d}}y{\mathrm{d}}\tau.
\end{eqnarray}
By Jensen's inequality, we have
\begin{eqnarray}\label{s2}
\int_{s}^{s+1}\!\!\int_{ B}\!|w(y,\tau
)|^{\frac{p+3}2}{\mathrm{d}}y{\mathrm{d}}\tau&\le&C\Big(\int_{s}^{s+1}\!\!\int_{
B}\!\!|w(y,\tau
)|^{p+1}{\mathrm{d}}y{\mathrm{d}}\tau\Big)^{\frac{p+3}{2(p+1)}}.
\end{eqnarray}
By combining  (\ref{s1}) and (\ref{s2}), we can write that
\begin{eqnarray}\label{s3}
\int_{ B}\!|w(y,\sigma (s))|^{\frac{p+3}2}{\mathrm{d}}y
&\le&C\Big(\int_{s}^{s+1}\!\!\int_{ B}\!\!|w(y,\tau
)|^{p+1}{\mathrm{d}}y{\mathrm{d}}\tau\Big)^{\frac{p+3}{2(p+1)}}.
\end{eqnarray}
Using   (\ref{s3}) and the fact that $ab\le a^2+b^2$, we have
\begin{eqnarray*}
\int_{ B}\!|w(y,s )|^{\frac{p+3}2}{\mathrm{d}}y&\le&\int_{
B}\!|w(y,\sigma (s))|^{\frac{p+3}2}{\mathrm{d}}y  +
C\int_{s}^{s+1}\int_{
B}\!|\partial_s w(y,\tau)||w(y,\tau)|^{\frac{p+1}2}{\mathrm{d}}y{\mathrm{d}}\tau\\
&\le& C\Big(\int_{s}^{s+1}\!\!\int_{ B}\!\!|w(y,\tau
)|^{p+1}{\mathrm{d}}y{\mathrm{d}}\tau\Big)^{\frac{p+3}{2(p+1)}}\\
&& + C\Big(\int_{s}^{s+1}\int_{
B}\!|w(y,\tau)|^{p+1}{\mathrm{d}}y{\mathrm{d}}\tau\Big)^{\frac12}
\Big(\int_{s}^{s+1}\int_{ B}\!|\partial_s
w(y,\tau)|^2{\mathrm{d}}y{\mathrm{d}}\tau\Big)^{\frac12}.
\end{eqnarray*}
Since $\ds{e^{-2\m s}\int_{s}^{s+1}\int_{ B}\!|w(y,\tau
)|^{p+1}{\mathrm{d}}y{\mathrm{d}}\tau\rightarrow 0}$ from
(\ref{F2}), we use (\ref{cor1}) to obtain (\ref{cor2}).

\medskip

\no - {\it Proof of \eqref{cor3}}: It follows from \eqref{cor2} through the Holder
inequality and (\ref{cor2}).

\medskip

\no - {\it Proof of \eqref{03bis}}:
Note that from now on, we assume that  $x_0$ is a non characteristic point with slope $\delta_0\in (0,1)$.
It is a direct consequence of (\ref{cor3})
and (\ref{cor00001}).
\medskip

\no - {\it Proof of \eqref{03}}:
Let $s\ge s_0+1$, $s_1=s_1(s)\in [s-1,s]$  and $s_2=s_2(s)\in
[s,s+1]$ to be chosen later. By integrating after multiplication by $e^{-2\m s}$ the
expression (\ref{e8}) of $I(w)$ in time between $s_1$ and $s_2$, we
obtain
\begin{eqnarray}\label{q2}
 &&\m\int_{s_1}^{s_2}e^{-2\m s}\int_{B}|\grad w|^2(1-|y|^2){\mathrm{d}}y{\mathrm{d}}s=
\underbrace{e^{-2\m s_2}I(w(s_2))-e^{-2\m s_1}I(w(s_1))}_{B_1(s)}\qquad \qquad \nonumber\\
&&+\underbrace{2\m \int_{s_1}^{s_2}e^{-2\m s}\int_{B} (y.\grad w)
\partial_sw {\mathrm{d}}y{\mathrm{d}}s}_{B_2(s)}
 +\underbrace{\eta \int_{s_1}^{s_2}e^{-2\m s}\int_{B}(\partial_s
w)^2{\mathrm{d}}y{\mathrm{d}}s}_{B_3(s)}\\
&&+ \underbrace{\m \int_{s_1}^{s_2}e^{-2\m
s}\int_{B}|w|^{p+1}{\mathrm{d}}y{\mathrm{d}}s}_{B_4(s)} \underbrace{
-(\m^2N  +\m \frac{2p+2}{(p-1)^2})\int_{s_1}^{s_2}e^{-2\m s}\int_{B}w^2{\mathrm{d}}y{\mathrm{d}}s}_{B_5(s)}\nonumber\\
&&\underbrace{- \m (\frac{p+3}{p-1}-N)\int_{s_1}^{s_2} e^{-2\m
s}\int_{B}w\partial_sw{\mathrm{d}}y{\mathrm{d}}s}_{B_6(s)}\underbrace{-
\int_{s_1}^{s_2}e^{-2\m s}\int_{\partial B}  (\partial_sw)^2
{\mathrm{d}}\sigma{\mathrm{d}}s}_{B_7(s)}\nonumber\\  &&+
\underbrace{\int_{s_1}^{s_2}e^{-2\m s}\int_{\partial B}
(\partial_sw-\m w)^2{\mathrm{d}}\sigma
{\mathrm{d}}s}_{B_8(s)}\underbrace{-\m\int_{s_1}^{s_2}e^{-2\m
s}\int_{B}(|y|^2|\grad w|^2-(y.\grad
w)^2){\mathrm{d}}y{\mathrm{d}}s}_{B_{9}(s)}.\nonumber
\end{eqnarray}
Now, we    control all the terms on the right-hand
 side of the
relation (\ref{q2}):

\no Note that, by (\ref{E2}) and using the Cauchy-Schwarz
inequality, we can write
\begin{equation}\label{I1}
e^{-2\m s_2}I(w(s_2))\le Ce^{-2\m s_2}
\displaystyle\int_{B}(\partial_s w(s_2))^2 {\mathrm{d}}y+Ce^{-2\m
s_2} \displaystyle\int_{B}w^2(s_2)
 {\mathrm{d}}y.
\end{equation}
By exploiting  (\ref{cor3}) and the fact that  $s_2\in [s,s+1]$, we
conclude that
\begin{equation}\label{I11}
e^{-2\m s_2} \displaystyle\int_{B}w^2(s_2)
 {\mathrm{d}}y\rightarrow 0\ \ {\textrm {as}}\ \ s\rightarrow
 +\infty,
\end{equation}
on the one hand. On the other hand, by   using the mean value
theorem, let us choose  $s_2=s_2(s)\in [s,s+1]$ such that
\begin{eqnarray}\label{q4}\displaystyle\int_{s}^{s+1}
e^{-2\m \tau}\displaystyle\int_{B}(\partial_s w(\tau))^2
{\mathrm{d}}y{\mathrm{d}}\tau&=&e^{-2\m
s_2}\displaystyle\int_{B}(\partial_s w(s_2))^2 {\mathrm{d}}y.
\end{eqnarray}
By combining  (\ref{03bis}) and (\ref{q4}) we
obtain
\begin{equation}\label{I111}
e^{-2\m s_2}\displaystyle\int_{B}(\partial_s w(s_2))^2
{\mathrm{d}}y\rightarrow 0\ \ {\textrm {as}}\ \ s\rightarrow
 +\infty.
\end{equation}
Then, by using (\ref{I1}), (\ref{I11}) and (\ref{I111}), we get
\begin{equation}\label{I1111}
e^{-2\m s_2}I(w(s_2))\rightarrow 0\ \ {\textrm {as}}\ \ s\rightarrow
+\infty.
\end{equation}
From (\ref{E2}) and the fact that $ab\le a^2+b^2$, we write
\begin{equation}\label{I02}
-e^{-2\m s_1}I(w(s_1))\le C e^{-2\m s_1}\displaystyle\int_{B}(\partial_s w(s_1))^2
{\mathrm{d}}y.
\end{equation}
Similarly, by   using the mean value
theorem, we choose $s_1=s_1(s)\in [s-1,s]$ such that
\begin{eqnarray}\label{q44}\displaystyle\int_{s-1}^{s}
e^{-2\m \tau}\displaystyle\int_{B}(\partial_s w(\tau))^2
{\mathrm{d}}y{\mathrm{d}}\tau&=&e^{-2\m
s_1}\displaystyle\int_{B}(\partial_s w(s_1))^2 {\mathrm{d}}y.
\end{eqnarray}
By  (\ref{03bis}), (\ref{I02}) and (\ref{q44}) we obtain
\begin{equation}\label{I2}
-e^{-2\m s_1}I(w(s_1))\rightarrow 0\ \ {\textrm {as}}\ \ s\rightarrow
+\infty.
\end{equation}
Note that by combining   (\ref{I1111}), (\ref{I2}) and the fact that $B_1(s)=e^{-2\m
s_2}I(w(s_2))-e^{-2\m s_1}I(w(s_1))$,
we deduce that
\begin{equation}\label{B1}
B_1(s)\rightarrow 0\ \ {\textrm {as}}\ \ s\rightarrow +\infty.
\end{equation}
To estimate $B_2(s)$, since $s_1\in [s-1,s]$ and $s_2\in [s,s+1]$, we write
\begin{eqnarray}\label{q01}
B_2(s)&\le&C\Big(\int_{s-1}^{s+1}e^{-2\m \tau}\int_{B} |\grad w|^2
{\mathrm{d}}y{\mathrm{d}}\tau \Big)^{\frac12}
\Big(\int_{s-1}^{s+1}e^{-2\m \tau}\int_{B} (
\partial_sw)^2 {\mathrm{d}}y{\mathrm{d}}\tau \Big)^{\frac12}.\qquad
\end{eqnarray}
By using (\ref{cor1}) and the covering argument of \cite{MZimrn05},
we have
\begin{eqnarray}\label{cor01}
e^{-2\eta s}\int_{s}^{s+1} \!\int_{B}\!|\grad
w|^2{\mathrm{d}}y{\mathrm{d}}\tau\le K.
\end{eqnarray}
Thus
\begin{eqnarray}\label{B02}
\int_{s-1}^{s+1}e^{-2\m \tau}\int_{B} |\grad w|^2
{\mathrm{d}}y{\mathrm{d}}\tau &\le C K.
\end{eqnarray}
Then, by (\ref{q01}), (\ref{B02}) and  (\ref{03bis}), we deduce
\begin{equation}\label{B2}
B_2(s)\rightarrow 0\ \ {\textrm {as}}\ \ s\rightarrow +\infty.
\end{equation}
By (\ref{03bis}), we can say that
\begin{equation}\label{B3}
B_3(s)\rightarrow 0\ \ {\textrm {as}}\ \ s\rightarrow +\infty.
\end{equation}
By (\ref{F2}), we also deduce that
\begin{equation}\label{B4}
B_4(s)\rightarrow 0\ \ {\textrm {as}}\ \ s\rightarrow +\infty.
\end{equation}
The terms $B_5(s)$ and $B_7(s)$  are negative. By
using  (\ref{cor3}) and (\ref{03bis}), we have
\begin{equation}\label{B6}
B_6(s)\rightarrow 0\ \ {\textrm {as}}\ \ s\rightarrow +\infty.
\end{equation}
By (\ref{F3}), we write that
\begin{equation}\label{B8}
B_8(s)\rightarrow 0\ \ {\textrm {as}}\ \ s\rightarrow +\infty.
\end{equation}
Finally, since $|y.\grad
w|\le |y| |\grad w|$, we can say that the term $B_9(s)$ is negative.
By combining (\ref{B1}), (\ref{B2}), (\ref{B3}), (\ref{B4}), (\ref{B6}), (\ref{B8})
and the fact that the terms $B_5(s)$, $B_7(s)$ and $B_9(s)$  are negative, we conclude that
\begin{eqnarray}\label{q9}
 \int_{s}^{s+1}e^{-2\m \tau}\int_{B}|\grad w|^2(1-|y|^2){\mathrm{d}}y{\mathrm{d}}\tau\rightarrow 0
 \ \ {\textrm {as}}\ \ s\rightarrow +\infty.
\end{eqnarray}
By using (\ref{q9}) and the covering argument of \cite{MZimrn05},
we deduce that  estimate (\ref{03})
holds.


\medskip

\no - {\it Proof of \eqref{04}}:
By integrating the functional $F(w,s)$ defined in
(\ref{E3}) in time between $s$ and $s+1$, we  write
\begin{eqnarray}\label{fin2}
\int_{s}^{s+1}\!\!F(w,\tau)
d\tau&=&\displaystyle\int_{s}^{s+1}\int_{B}\!\!e^{-2\eta \tau}\Big (
\frac{1}{2}(\partial_s w)^2
+\frac{p+1}{(p-1)^2}w^2-\frac{1}{p+1}|w|^{p+1}\Big ) {\mathrm{d}}y{\mathrm{d}}\tau\nonumber\\
&&+\frac{1}{2}\displaystyle\int_{s}^{s+1}\!\!e^{-2\eta
\tau}\int_{B}\!\!\Big (|\grad w|^2-(y.\grad w)^2\Big )
{\mathrm{d}}y{\mathrm{d}}\tau\\
&&-\eta\int_{s}^{s+1}\!\!e^{-2\eta
\tau}\displaystyle\int_{B}w\partial_s w
{\mathrm{d}}y{\mathrm{d}}\tau+ \frac{\m N}2
\int_{s}^{s+1}\!\!e^{-2\eta \tau}\displaystyle\int_{B}w^2
 {\mathrm{d}}y{\mathrm{d}}\tau.\nonumber
\end{eqnarray}
By using (\ref{cor3}),  (\ref{03}), (\ref{F2}) and  (\ref{fin2}), we conclude that
\begin{eqnarray}\label{fin3}
\int_{s}^{s+1}\!\!F(w,\tau)
d\tau\rightarrow 0,
 \ \ {\textrm {as}}\ \ s\rightarrow +\infty.
\end{eqnarray}
Combining The monotonicity of $F(w,s)$ proved in Proposition \ref{pro}, and (\ref{fin3}), we deduce the identity (\ref{04}).  This concludes
 the proof
 of  Proposition \ref{final}.

\Box


\no Since the derivation of Theorem \ref{t0}   from Proposition  \ref{final}  is the
same as in  \cite{MZimrn05} (up to some very minor changes),
this concludes the proof of Theorem \ref{t0}.
\Box

\def\cprime{$'$}

\noindent{\bf Address}:\\
Universit\'e de Tunis El-Manar, Facult\'e des Sciences de Tunis,
D\'epartement de math\'ematiques, Campus Universitaire 1060,
 Tunis, Tunisia.\\
\vspace{-7mm}
\begin{verbatim}
e-mail: ma.hamza@fst.rnu.tn
\end{verbatim}
Universit\'e Paris 13, Sorbonne Paris Cit\'e, LAGA, CNRS (UMR 7539),
99 avenue J.B. Cl\'ement, 93430 Villetaneuse, France.\\
\vspace{-7mm}
\begin{verbatim}
e-mail: Hatem.Zaag@univ-paris13.fr
\end{verbatim}
\end{document}